\documentclass[a4paper,11pt]{article}
\usepackage[top=2.5cm,bottom=2.5cm,left=2.2cm,right=2.2cm]{geometry}
\usepackage{amsfonts}
\usepackage{mathrsfs,amscd,amssymb,amsthm,amsmath,bm,graphicx,psfrag,subfigure,url,mathtools}
\usepackage{pict2e}
\usepackage{stfloats}
\usepackage{psfrag,amsmath}
\usepackage{tikz}
\usepackage{indentfirst}
\usepackage{hyperref}
\usepackage{bookmark}
\usepackage{enumerate}
\usepackage{latexsym,euscript,epic,eepic,color}
\usepackage{multirow}
\usepackage{multicol}
\usepackage{longtable}
\usepackage{adjustbox}
\usepackage[all]{xy}
\usepackage{setspace}
\usepackage{epstopdf}
\allowdisplaybreaks
\usepackage{authblk}
\usepackage{tikz}
\pgfdeclarelayer{nodelayer}
\pgfdeclarelayer{edgelayer}
\pgfsetlayers{nodelayer,edgelayer,main}

\makeatletter

\renewcommand{\@seccntformat}[1]{{\csname the#1\endcsname}{\normalsize .}\hspace{.5em}}
\makeatother

\usepackage{ifpdf}
\usepackage{indentfirst}

\def \[{\begin{equation}}
\def \]{\end{equation}}

\newtheorem{theorem}{Theorem}

\newtheorem{claim}{Claim}
\newtheorem{lemma}{Lemma}
\newtheorem{corollary}{Corollary}

\newtheorem{problem}{Problem}

\newtheorem{conjecture}{Conjecture}
\newtheorem{remark}{Remark}

\newenvironment{wst}
{\setlength{\leftmargini}{1.5\parindent}
 \begin{itemize}
 \setlength{\itemsep}{-1.1mm}}
{\end{itemize}}
\begin{document}
\baselineskip=0.23in
\title{\bf Tur\'an-type problems on $[a,b]$-factors\\ of graphs, and beyond\thanks{S.L. financially supported by the National Natural Science Foundation of China (Grant Nos. 12171190, 11671164) and the Special Fund for Basic Scientiﬁc Research of Central Colleges (Grant No. CCNU24JC005).\\[3pt] \hspace*{5mm}{\it Email addresses}: yfhaomath@sina.com (Y. Hao),\ li@ccnu.edu.cn (S. Li)}}
\author[1]{Yifang Hao}
\author[,1,2]{Shuchao Li\thanks{Corresponding author}}
\affil[1]{School of Mathematics and Statistics, and Hubei Key Lab--Math. Sci.,\linebreak Central China Normal University, Wuhan 430079, China}
\affil[2]{Key Laboratory of Nonlinear Analysis \& Applications (Ministry of Education),\linebreak Central China Normal University, Wuhan 430079, China}
\date{\today}
\maketitle
\begin{abstract}
 Given a set of graphs $\mathcal{H}$, we say that a graph $G$ is \textit{$\mathcal{H}$-free} if it does not contain any member of $\mathcal{H}$ as a subgraph. Let $\text{ex}(n,\mathcal{H})$ (resp. $\text{ex}_{sp}(n,\mathcal{H})$) denote the maximum size (resp. spectral radius) of an $n$-vertex $\mathcal{H}$-free graph. Denote by $\text{Ex}(n, \mathcal{H})$ the set of all $n$-vertex $\mathcal{H}$-free graphs with $\text{ex}(n, \mathcal{H})$ edges. Similarly,  let $\mathrm{Ex}_{sp}(n,\mathcal{H})$ be the set of all $n$-vertex $\mathcal{H}$-free graphs with spectral radius $\text{ex}_{sp}(n, \mathcal{H})$. For positive integers $a, b$ with $a\leqslant b$, an $[a,b]$-factor of a graph $G$ is a spanning subgraph $F$ of $G$ such that $a\leqslant d_F(v)\leqslant b$ for all $v\in V(G)$, where $d_F(v)$ denotes the degree of the vertex $v$ in $F.$ Let $\mathcal{F}_{a,b}$ be the set of all the $[a,b]$-factors of an $n$-vertex complete graph $K_n$. In this paper, we determine the Tur\'an number $\text{ex}(n,\mathcal{F}_{a,b})$ and the spectral Tur\'an number $\text{ex}_{sp}(n,\mathcal{F}_{a,b}),$ respectively. Furthermore, the bipartite analogue of $\text{ex}(n,\mathcal{F}_{a,b})$ (resp.  $\text{ex}_{sp}(n,\mathcal{F}_{a,b})$) is also obtained. All the corresponding extremal graphs are identified. Consequently, one sees that $\mathrm{Ex}_{sp}(n,\mathcal{F}_{a,b})\subseteq \text{Ex}(n, \mathcal{F}_{a,b})$ holds for graphs and bipartite graphs. This partially answers an open problem proposed by Liu and Ning \cite{LN2023}. Our results may deduce a main result of Fan and Lin \cite{FL2022}.\\

\vspace{0.2cm}

\end{abstract}

\section{\normalsize Introduction}

Extremal graph theory, one of the most important branches in combinatorics, intends to study how global properties of a graph control its local structure. Tur\'an-type problem is a typical representative in extremal graph theory. Given a set of graphs $\mathcal{H}$, a graph $G$ is \textit{$\mathcal{H}$-free} if it contains no member of $\mathcal{H}$ as a subgraph. In particular, if $\mathcal{H}=\{H\}$, then we also say that $G$ is $H$-free. Let $\text{ex}(n,\mathcal{H})$  denote the maximum size of an $n$-vertex $\mathcal{H}$-free graph. Denote by $\text{Ex}(n, \mathcal{H})$ the set of all $n$-vertex $\mathcal{H}$-free graphs with $\text{ex}(n, \mathcal{H})$ edges. One of the central problems in extremal graph theory is to study the behavior of the Tur\'an number $\text{ex}(n,\mathcal{H})$ and to characterize all the graphs in $\text{Ex}(n, \mathcal{H})$. The classical results in this area include Mantel's theorem \cite{WM1907} that any $n$-vertex graph with more than $\lfloor\frac{n^2}{4}\rfloor$ edges must contain a triangle. The cornerstone result in extremal graph theory is Tur\'an's theorem \cite{PT1941} in which Tur\'an determined $\text{ex}(n, K_r)$ in 1941. Five years later, the celebrated  Erd\H{o}s-Stone theorem \cite{ES1966,ES1946} presented an asymptotic solution for $\text{ex}(n, H)$ when $\chi(H)\geqslant 3$. For more development along this line one may consult the nice paper \cite{FS2013}.

Spectral extremal graph theory, comparing with the classical extremal graph theory, is much younger. In the past thirty years, it has experienced rapid development. Given a set of  graphs $\mathcal{H}$, let $\text{ex}_{sp}(n,\mathcal{H})$ denote the maximum adjacency spectral radius of an $n$-vertex $\mathcal{H}$-free graph, which is the so-called \textit{spectral Tur\'an number}. Denote by $\mathrm{Ex}_{sp}(n,\mathcal{H})$ the set of all $n$-vertex $\mathcal{H}$-free graphs with adjacency spectral radius $\text{ex}_{sp}(n, \mathcal{H})$. Nikiforov \cite{NV2011-1} initiated the research on the spectral Tur\'an-type problems: Determine the maximum adjacency spectral radius over the class of $n$-vertex $\mathcal{H}$-free graphs and characterize the corresponding extremal graphs. Subsequently, he conducted a systematic research on this topic, and spectral extremal graph theory attracts more and more researchers' attention from now on. One may consult the nice survey \cite{LFL2022} for more details.

Clearly, the spectral Tur\'an-type problem has close relationship with Tur\'an-type problem. Both of them have the same goal: to determine both $\text{ex}(n,\mathcal{H})$ and $\text{ex}_{sp}(n,\mathcal{H})$ and identify the extremal graphs in $\mathrm{Ex}(n,\mathcal{H})$ and $\mathrm{Ex}_{sp}(n,\mathcal{H})$, respectively. One sees, from the mathematical literature, the achievement on Tur\'an-type problems is much richer than that of spectral Tur\'an-type problems. So it is very natural for us to deduce the spectral analogues of the Tur\'an-type problems. On the other hand, some mathematical phenomenons reveal that $\mathrm{Ex}_{sp}(n,\mathcal{H})$ has close relation with $\mathrm{Ex}(n,\mathcal{H})$. This leads us to reveal the mysterious veil between $\mathrm{Ex}(n,\mathcal{H})$ and $\mathrm{Ex}_{sp}(n,\mathcal{H})$.

We will survey some typical results surrounding the above observation in what follows. Nikiforov \cite{Nik2007} and Guiduli \cite{Gui1996}, independently, determined $\text{ex}_{sp}(n, K_{r+1})$. Together with \cite{PT1941}, one sees $\mathrm{Ex}_{sp}(n,K_{r+1})\subseteq \mathrm{Ex}(n,K_{r+1})$. Erd\H{o}s, F\"{u}redi, Gould, and Gunderson \cite{Erdo} determined the Tur\'an number $\mathrm{ex}(n, K_1\vee kK_2)$ and identified the extremal graphs, whose spectral analogue is obtained in \cite{Cioa,Zhai}, in which Cioab\u{a}, Feng, Tait and Zhang \cite{Cioa} showed that $\mathrm{Ex}_{sp}(n,K_1\vee kK_2)\subseteq \mathrm{Ex}(n,K_1\vee kK_2)$ and Zhai, Liu and Xue \cite{Zhai} determined the unique extremal graph in $\mathrm{Ex}_{sp}(n,K_1\vee kK_2)$. Chen, Gould, Pfender and Wei \cite{Cheng} generalized the main result in \cite{Erdo} determining the Tur\'an number $\mathrm{ex}(n, K_1\vee kK_r)$ and identifying the extremal graphs. Its spectral analogue is obtained in \cite{Desai,You}, in which You, Wang and Kang \cite{You} determined the unique extremal graph in $\mathrm{Ex}_{sp}(n,K_1\vee kK_r)$ and Desai et al. \cite{Desai} showed $\mathrm{Ex}_{sp}(n, K_1\vee kK_r)\subseteq \mathrm{Ex}(n, K_1\vee kK_r)$ for sufficiently large $n$. Let $F_1, \ldots, F_t$ be $t$ disjoint color-critical graphs with $\chi(F_i)=r+1\, (r\geqslant 2)$. Simonovits \cite{Simo} determined the Tur\'an number $\mathrm{ex}(n, \cup_{i=1}^tF_i)$ and identified the extremal graph for sufficiently large $n$. Recently, Lei and Li \cite{Lei} determined the spectral Tur\'an number $\mathrm{ex}_{sp}(n, \cup_{i=1}^tF_i)$ and identified the extremal graph for sufficiently large $n$. One sees that the unique extremal graph in $\mathrm{Ex}_{sp}(n, \cup_{i=1}^tF_i)$ coincides with that of $\mathrm{Ex}(n, \cup_{i=1}^tF_i)$. For more advances on this topic, we refer the reader to \cite{BDT2024,HLL2024,SLW2023,SLW2023-1,Wang2,Zhai-b} and the nice survey paper \cite{LFL2022}.

Based on the above achievements and some other mathematical phenomenon, Liu and Ning \cite{LN2023} proposed an interesting problem as follows.
\begin{problem}[\cite{LN2023}]\label{pb1}
Let $H$ be any graph. Characterize all graphs $H$ such that
$$
\text{$\mathrm{Ex}_{sp}(n,H)\subseteq \mathrm{Ex}(n,H)$}
$$
for sufficiently large $n$.
\end{problem}

Recently, there are two breakthroughs for Problem \ref{pb1}. The first one is due to Wang, Kang and Xue's work \cite{Wang}: Given a graph $H$ with $\mathrm{ex}(n, H)=t_r(n) +O(1)$, where $r \geqslant 2$ and $t_r(n)$ denotes the size of Tur\'an graph $T_r(n)$, then $\mathrm{Ex}_{sp}(n,H)\subseteq \mathrm{Ex}(n,H)$ when $n\rightarrow \infty.$ This confirms a conjecture proposed by Cioab\u{a}, Desai, and Tait \cite[Conjecture 7.1]{Cioa2}. Another one is due to Byrne, Desai and Tait's work \cite{BDT2024}. It proves a general theorem to characterize the spectral extremal graphs for a wide range of forbidden families $\mathcal{H}$ and implies several new and existing results. Particularly, \cite{BDT2024} deduces the following: Whenever $\mathrm{ex}(n, \mathcal{H})=O(n)$, $K_{k+1,\infty}$ is not $\mathcal{H}$-free and $\mathrm{Ex}(n, \mathcal{H})$ contains the complete bipartite graph
$K_{k,n-k}$ (or certain similar graphs), then $\mathrm{Ex}_{sp}(n, \mathcal{H})\subseteq\mathrm{Ex}(n, \mathcal{H})$ for sufficiently large $n$.

In order to establish a criterion for $\mathrm{Ex}_{sp}(n,H)\cap \mathrm{Ex}(n,H)\not=\emptyset,$ they gave a conjecture:
\begin{conjecture}\label{conj-1}
Let $k$ be a fixed positive integer and $n$ be a sufficiently large integer. Let $F$
be a graph such that $\text{\rm ex}(n,F)=\frac{1}{2}n^2-kn+O(1).$ Then we have $\mathrm{Ex}_{sp}(n,F)\subseteq \mathrm{Ex}(n,F)$.
\end{conjecture}

In this paper, motivated by \cite{LM-2021,L-10,ML2023,LN2023}, we contribute to Problem \ref{pb1} by proving positive results when $\mathcal{F}_{a,b}$ is the set of all the $[a,b]$-factors of a complete graph $K_n$, and we also contribute to Problem \ref{pb1} by proving positive results when $\mathcal{B}_{a,b}$ is the set of all the $[a,b]$-factors of bipartite graphs on $n$ vertices. Our results also present certain relation with Conjecture~\ref{conj-1}.

Our first main result determines the maximum size of an $n$-vertex graph forbidding $[a,b]$-factors, and characterizes the extremal graphs.
\begin{theorem}\label{thm1.1}
Let $a\leqslant b$ be two positive integers, and $G$ be a graph of order $n$, where $n\geqslant a+1$ and $na\equiv0\pmod{2}$ when $a=b$. If $G$ contains no  $[a,b]$-factors, then $e(G)\leqslant\binom{n-1}{2}+a-1$ with equality if and only if one of the following holds:
\begin{wst}
\item[{\rm (i)}]$G\cong K_{a-1}\vee (K_{n-a}\cup K_1)$ or $K_{1,3},$ if $ab=1$ or $ab=2;$
\item[{\rm (ii)}]$G\cong K_{a-1}\vee (K_{n-a}\cup K_1)$ or $K_2\vee \overline{K_3},$ if $a=b=2;$
\item[{\rm (iii)}]$G\cong K_{a-1}\vee (K_{n-a}\cup K_1),$ if $b\geqslant3$.
\end{wst}
\end{theorem}

Our second main result determines the maximum adjacency spectral radius of an $n$-vertex graph forbidding $[a,b]$-factors and characterizes the extremal graphs, which strengthens the main result of Wei and Zhang \cite[Theorem 1]{WZ2023}.
\begin{theorem}\label{thm1.2}
Let $a\leqslant b$ be two positive integers, and $G$ be a graph of order $n$, where $n\geqslant a+1$ and $na\equiv0\pmod{2}$ when $a=b$. If $G$ contains no $[a,b]$-factors, then
$
\rho(G)\leqslant \rho(K_{a-1}\vee (K_{n-a}\cup K_1))
$
with equality if and only if $G\cong K_{a-1}\vee (K_{n-a}\cup K_1)$.
\end{theorem}

An immediate consequence of Theorems~\ref{thm1.1} and \ref{thm1.2} directly contributes to Problem~\ref{pb1}.
\begin{corollary}
Let $a,b,n$ be three positive integers with $a\leqslant b$, $n\geqslant a+1$ and $na\equiv0\pmod{2}$ when $a=b$. Then $\mathrm{Ex}_{sp}(n, \mathcal{F}_{a,b})\subseteq\mathrm{Ex}(n, \mathcal{F}_{a,b})$.
\end{corollary}

In what follows, we give the bipartite analogues of Theorems \ref{thm1.1} and \ref{thm1.2}. In fact, they are closely relative to the \textit{Zarankiewicz problem}, which asks how many edges an $n$ by $m$ bipartite graph may have without containing a copy of $K_{s,t}$, and is the other most famous bipartite Tur\'an problem.

Before formulating our main results, we recall the definition of double nested graph (see \cite{A.M.,ML2023,LS2020,P.M.}), which is also called the \textit{bipartite chain graph} \cite{BFP2008}.
\begin{figure}
\begin{center}
\includegraphics[width=50mm]{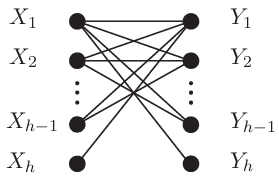} \\
  \caption{The structure of a double nested graph.}\label{fig1}
\end{center}
\end{figure}

Let $G=(X,Y)$ be a connected bipartite graph. We call $G$ a \textit{double nested graph} if there exist partitions $X=X_1\cup X_2\cup \cdots\cup X_h$ and $Y=Y_1\cup Y_2\cup \cdots\cup Y_h$ such that all vertices in $X_i$ are adjacent to all vertices in $\bigcup_{j=1}^{h+1-i}Y_j$ for $1\leqslant i\leqslant h$ (see Fig.~\ref{fig1}), in which each solid circle denotes an independent set of an appropriate size, each line between two large solid circles means that all vertices in one large solid circle are adjacent to all vertices in the other one. For convenience, let $|X_i|=p_i$ and $|Y_i|=q_i$ for $i=1,2,\ldots,h.$ Then we denote the graph $G$ by
$
D(p_1,p_2,\ldots,p_h;q_1,q_2,\ldots,q_h).
$
Clearly, if $p_1=0$ or $q_1=0$, then $
D(p_1,p_2,\ldots,p_h;q_1,q_2,\ldots,q_h).
$ is the disjoint union of a connected graph and some isolated vertices. Hence, a double nested graph $D(p_1,p_2,\ldots,p_h;q_1,q_2,\ldots,q_h)$ is connected if and only if $p_1,q_1>0$.

The next main result determines the maximum size of an $|X|$ by $|Y|$ bipartite graph without containing any $[a,b]$-factor, and the extremal graphs are also characterized.
\begin{theorem}\label{thm1.3}
Let $a\leqslant b$ be two positive integers. Let $G=(X,Y)$ be a bipartite graph with $|X|\leqslant|Y|$ forbidding $[a,b]$-factors.
\begin{wst}
\item[{\rm (i)}]If $a|Y|>b|X|$, then $e(G)\leqslant |X||Y|$ with equality if and only if $G\cong K_{|X|,|Y|}$.
\item[{\rm (ii)}]If $a|Y|\leqslant b|X|$ and $a>|X|$, then
$e(G)\leqslant |X||Y|$ with equality if and only if $G\cong K_{|X|,|Y|}$.
\item[{\rm (iii)}]If $a|Y|\leqslant b|X|$ and $a\leqslant |X|$, then
$
e(G)\leqslant |X|(|Y|-1)+a-1
$
with equality if and only if $G\cong D(a-1,|X|-a+1;|Y|-1,1)$.
\end{wst}
\end{theorem}

The next main result determines the maximum size of an $n$-vertex bipartite graph without containing any $[a,b]$-factor, and characterizes the corresponding extremal graphs. For convenience, let 
\begin{equation}\label{eq:001}
f(a,b):=\left\lfloor\frac{an-1}{a+b}\right\rfloor\left(n-\left\lfloor\frac{an-1}{a+b}\right\rfloor\right).
\end{equation}
\begin{theorem}\label{thm1.4}
Let $a,b,n$ be three positive integers with $a\leqslant b$ and $a\leqslant\lfloor\frac{n}{2}\rfloor$. Let $G$ be an $n$-vertex bipartite graph forbidding $[a,b]$-factors. 
\begin{wst}
\item[{\rm (i)}]If $f(a,b)>\left\lfloor\frac{n}{2}\right\rfloor(\left\lceil\frac{n}{2}\right\rceil-1)+a-1$, then 
    $
    e(G)\leqslant f(a,b)
    $
    with equality if and only if $G\cong K_{\lfloor\frac{an-1}{a+b}\rfloor,n-\lfloor\frac{an-1}{a+b}\rfloor}$.
\item[{\rm (ii)}]If  $f(a,b)=\frac{n}{2}(\frac{n}{2}-1)+a-1$ for even $n$, then 
    $
    e(G)\leqslant f(a,b)
    $
    with equality if and only if $G\in\{K_{\lfloor\frac{an-1}{a+b}\rfloor,n-\lfloor\frac{an-1}{a+b}\rfloor}, D(a-1,\frac{n}{2}-a;\frac{n}{2},1), D(a-1,\frac{n}{2}-a+1;\frac{n}{2}-1,1)\}$.
\item[{\rm (iii)}]If  $f(a,b)=(\frac{n-1}{2})^2+a-1$ for odd $n$, then 
    $
    e(G)\leqslant f(a,b)
    $
    with equality if and only if $G\in\{K_{\lfloor\frac{an-1}{a+b}\rfloor,n-\lfloor\frac{an-1}{a+b}\rfloor}, D(a-1,\frac{n+1}{2}-a;\frac{n-1}{2},1)\}$.
\item[{\rm (iv)}]If  $f(a,b)<\frac{n}{2}(\frac{n}{2}-1)+a-1$ for even $n$, then $
    e(G)\leqslant \frac{n}{2}(\frac{n}{2}-1)+a-1$ 
    with equality if and only if $G\in\{D(a-1,\frac{n}{2}-a;\frac{n}{2},1), D(a-1,\frac{n}{2}-a+1;\frac{n}{2}-1,1)\}$.
\item[{\rm (v)}]If  $f(a,b)<(\frac{n-1}{2})^2+a-1$ for odd $n$, then 
    $
    e(G)\leqslant (\frac{n-1}{2})^2+a-1
    $
    with equality if and only if $G\cong D(a-1,\frac{n+1}{2}-a;\frac{n-1}{2},1)$.
\end{wst}
\end{theorem}

The following main result determines the maximum adjacency spectral radius of an $|X|$ by $|Y|$ bipartite graph without containing any $[a,b]$-factor, and identifies the corresponding extremal graphs.
\begin{theorem}\label{thm1.5}
Let $a\leqslant b$ be two positive integers. Let $G=(X,Y)$ be a bipartite graph with $|X|\leqslant|Y|$ forbidding $[a,b]$-factors.
\begin{wst}
\item[{\rm (i)}]If $a|Y|>b|X|,$ then $\rho(G)\leqslant \sqrt{|X||Y|}$ with equality if and only if $G\cong K_{|X|,|Y|}.$
\item[{\rm (ii)}]If $a|Y|\leqslant b|X|$ and $a>|X|,$ then
$\rho(G)\leqslant \sqrt{|X||Y|}$ with equality if and only if $G\cong K_{|X|,|Y|}$.
\item[{\rm (iii)}]If $a|Y|\leqslant b|X|$ and $a\leqslant |X|,$ then
$
\rho(G)\leqslant\rho(D(a-1,|X|-a+1;|Y|-1,1))
$
with equality if and only if $G\cong D(a-1,|X|-a+1;|Y|-1,1)$.
\end{wst}
\end{theorem}

Let $G$ be a connected balanced bipartite graph of order $n$ and let $2\leqslant a=b=k\leqslant \frac{n}{2}-1$ in Theorem~\ref{thm1.5}. Consequently, our result deduces a main result of Fan and Lin \cite[Theorem 1.3]{FL2022}, which provides a sufficient spectral condition for a connected balanced bipartite graph to contain a $k$-factor.
\begin{theorem}[\cite{FL2022}]
Let $2\leqslant k\leqslant \frac{n}{2}-1$ and let $G$ be a connected balanced bipartite graph of order $n$. If
$$
\rho(G)\geqslant \rho(D(k-1,\frac{n}{2}-k+1;\frac{n}{2}-1,1)),
$$
then $G$ has a $k$-factor, unless $G\cong D(k-1,\frac{n}{2}-k+1;\frac{n}{2}-1,1)$.
\end{theorem}

Our last main result determines the maximum spectral radius of an $n$-vertex bipartite graph without $[a,b]$-factors, and identifies the corresponding extremal graphs. Recall that $f(a,b)$ is given in \eqref{eq:001}.
\begin{theorem}\label{thm1.7}
Let $a,b,n$ be three positive integers with $a\leqslant b$ and $a\leqslant\lfloor\frac{n}{2}\rfloor$. Let $G$ be a bipartite graph of order $n$ forbidding $[a,b]$-factors. 
\begin{wst}
\item[{\rm (i)}]If $\rho(D(a-1,\lceil\frac{n}{2}\rceil-a;\lfloor\frac{n}{2}\rfloor,1))<\sqrt{f(a,b)},$ then 
    $
    \rho(G)\leqslant \sqrt{f(a,b)}
    $
    with equality if and only if $G\cong K_{\lfloor\frac{an-1}{a+b}\rfloor,n-\lfloor\frac{an-1}{a+b}\rfloor}$.
\item[{\rm (ii)}]If $\rho(D(a-1,\lceil\frac{n}{2}\rceil-a;\lfloor\frac{n}{2}\rfloor,1))=\sqrt{f(a,b)},$ then 
    $
    \rho(G)\leqslant \sqrt{f(a,b)}
    $
    with equality if and only if $G\cong K_{\lfloor\frac{an-1}{a+b}\rfloor,n-\lfloor\frac{an-1}{a+b}\rfloor}$ or $D(a-1,\lceil\frac{n}{2}\rceil-a;\lfloor\frac{n}{2}\rfloor,1)$.
\item[{\rm (iii)}]If $\rho(D(a-1,\lceil\frac{n}{2}\rceil-a;\lfloor\frac{n}{2}\rfloor,1))>\sqrt{f(a,b)},$ then
    $
    \rho(G)\leqslant\rho(D(a-1,\left\lceil\frac{n}{2}\right\rceil-a;\left\lfloor\frac{n}{2}\right\rfloor,1))
    $
    with equality if and only if
    $G\cong D(a-1,\lceil\frac{n}{2}\rceil-a;\lfloor\frac{n}{2}\rfloor,1)$.
\end{wst}
\end{theorem}
Let $\text{BEx}(n, \mathcal{H})$ be the set of all $n$-vertex $\mathcal{H}$-free bipartite graphs with maximum edges, and  let $\mathrm{BEx}_{sp}(n,\mathcal{H})$ be the set of all $n$-vertex $\mathcal{H}$-free bipartite graphs with maximum spectral radius.  An immediate consequence of Theorem~\ref{thm1.5} and Theorem~\ref{thm1.7} contributes to Problem~\ref{pb1}.
\begin{corollary}
Let $a,b,n$ be three positive integers with $a\leqslant b$ and $a\leqslant\lfloor\frac{n}{2}\rfloor$. Then $\mathrm{BEx}_{sp}(n, \mathcal{B}_{a,b})\subseteq\mathrm{BEx}(n, \mathcal{B}_{a,b})$.
\end{corollary}

\vspace{3mm}

\noindent{\bf Outline of the paper}\ \ In section 2, some necessary preliminaries are given.  In Section 3, we give the proofs of Theorems~\ref{thm1.1} and \ref{thm1.2}.
In fact, we determine $\mathrm{ex}(n, \mathcal{F}), \ \mathrm{ex_{sp}}(n, \mathcal{F})$ and corresponding $\mathrm{Ex}(n, \mathcal{F}), \ \mathrm{Ex_{sp}}(n, \mathcal{F}),$ where
$\mathcal{F}$ is the set of all the $[a,b]$-factors of an $n$-vertex graph $G$. In Section 4, we firstly give the proof of Theorems~\ref{thm1.3}. Then based on Theorem~\ref{thm1.3}, we give the proof of \ref{thm1.4}. In Section 5, we give the proof of Theorems~\ref{thm1.5} at first. Then based on Theorem~\ref{thm1.5}, we give the proof of Theorem~\ref{thm1.7}. Some concluding remarks are given in the last section.\\

\noindent{\bf Notations and definitions}\ \
In this paper, we consider only finite, simple and undirected graphs. For graph theoretic notation and terminology not defined here, we refer the reader to \cite{GR2001,D.B.}.

Let $G=(V(G), E(G))$ be a graph with vertex set $V(G)$ and edge set $E(G)$. The \textit{order} of $G$ is the number $n:=|V(G)|$ of its vertices, and the \textit{size} of $G$ is the number $e(G):=|E(G)|$ of its edges. For a vertex $v\in V(G)$, the \textit{degree} of $v$, denoted by $d_G(v)$, is the number of vertices adjacent to $v$ in $G.$ The \textit{minimum degree} $\delta(G)=\min\{d_G(v): v\in V(G)\}.$
A graph $G$ is \textit{Hamiltonian} if it contains a Hamilton cycle, i.e., a cycle containing all vertices of $G$.

For two graphs $G$ and $H$, we define $G\cup H$ to be their \textit{disjoint union}. We write $tG$ to denote the disjoint union of $t$ copies of $G$. The \textit{join} of $G$ and $H$, denoted by $G \vee H$, is the graph obtained from $G\cup H$ by adding edges joining every vertex of $G$ to every vertex of $H$. A graph is \textit{color-critical} if it contains an edge whose deletion reduces its chromatic number.

The \textit{adjacency matrix} $A(G)=(a_{ij})$ of $G$ is defined as an $n\times n$ $(0,1)$-matrix with $a_{ij}=1$ if and only if $ij\in E(G)$. Note that $A(G)$ is real symmetric, its eigenvalues $\lambda_i$ are real. So we can index them as $\lambda_1\geqslant \lambda_2\geqslant\cdots\geqslant \lambda_n.$ The largest eigenvalue $\lambda_1$ of $A(G)$ is called the \textit{spectral radius} of $G$, written as $\rho(G)$.

For two disjoint subsets $A, B\subseteq V(G),$ we use $E(A,B)$ to denote the set of edges with one endpoint in $A$ and the other endpoint in $B$, and let $e(A,B)=|E(A,B)|.$  
Given a graph $G$, let $g,f: V(G)\rightarrow \mathbb{Z}$ be two functions with $g(v)\leqslant f(v)$ for all $v\in V(G).$ A \textit{$(g,f)$-factor} of the graph $G$ is a spanning subgraph $F$ of $G$ such that $g(v)\leqslant d_F(v)\leqslant f(v)$ for all $v\in V(G)$.
Let $a\leqslant b$ be two positive integers. If $g(v)=a$ and $f(v)=b$ for all $v\in V(G),$ then a $(g,f)$-factor is called an \textit{$[a,b]$-factor}. In particular, if $a=b=k$, then a $[k,k]$-factor is called a \textit{$k$-factor}.

\section{\normalsize Preliminaries}
In this section, we give some preliminaries, which will be used in the subsequent sections.
Consider a real matrix $M$ whose rows and columns are indexed by $V=\{1, \ldots, n\}.$ Assume that $M$ can be written as
$$
M=\left(
  \begin{array}{ccc}
    M_{11} & \cdots & M_{1s} \\
    \vdots & \ddots & \vdots \\
    M_{s1} & \cdots & M_{ss} \\
  \end{array}
  \right)
$$
according to a vertex partition $\pi:\ V=V_1\cup\cdots\cup V_s,$ wherein $M_{ij}$ denotes the submatrix (block) of $M$ formed by rows in $V_i$ and columns in $V_j.$ The \textit{quotient matrix} of $M$ is the matrix whose entries are the average row sums of the blocks of $M.$ The partition is called \textit{equitable} if each block $M_{ij}$ has a constant row sum.

\begin{lemma}[\cite{A.E.,YLH}]\label{lem2.1}
Let $M$ be a real square matrix with an equitable partition $\pi,$ and let $M_{\pi}$ be the corresponding quotient matrix. Then every eigenvalue of $M_{\pi}$ is an eigenvalue of $M.$ Furthermore, if $M$ is nonnegative, then the largest eigenvalues of $M$ and $M_{\pi}$ are equal.
\end{lemma}

\begin{lemma}[\cite{R.B.}]\label{lem2.2}
Let $G$ be a connected graph, and $H$ be a subgraph of $G$. Then $\rho(H)\leqslant\rho(G)$ with equality if and only if $H\cong G$.
\end{lemma}

The following result is an immediate consequence of Lemma \ref{lem2.2}.
\begin{corollary}\label{cor2.3}
Let $G=G'\cup lK_1$ be the disjoint union of a connected graph $G'$ and $l$ isolated vertices, where $l\geqslant0$. If $H$ is a spanning subgraph of $G$, then $\rho(H)\leqslant\rho(G)$ with equality if and only if $H\cong G$.
\end{corollary}

The following lemma provides a relation between the spectral radius of a graph $G$ and the size of $G$, which plays an important role in our later proof.
\begin{lemma}[\cite{Y.H.}]\label{lem2.4}
Let $G$ be a connected graph of order $n$. Then $\rho(G)\leqslant\sqrt{2e(G)-n+1}$ with equality if and only if $G\cong K_{1,n-1}$ or $G\cong K_n$.
\end{lemma}

\begin{lemma}[\cite{R.P.}]\label{lem2.5}
Let $\mathbb{G}(n,e)$ be the set of graphs with $n$ vertices and $e$ edges. If $e=$
\begin{math}
\left( \begin{smallmatrix}
r\\2
\end{smallmatrix} \right)+t,
\end{math}
where $0<t\leqslant r$, then $(K_t\vee (K_{r-t}\cup K_1))\cup (n-r-1)K_1$ is the unique graph with maximum spectral radius among all graphs in $\mathbb{G}(n,e)$.
\end{lemma}

Let $p,q,e$ be positive integers with $p\leqslant q$ and $e\leqslant pq$. Let $\mathcal{K}(p,q,e)$ be the set of bipartite graphs $G=(X,Y)$ with $|X|=p$, $|Y|=q$ and $e$ edges. The following lemma gives an upper bound on the spectral radius of graphs in $\mathcal{K}(p,q,e)$ with restriction on the value of $e$.
\begin{lemma}[\cite{LW2015}]\label{lem2.6}
If $p,q,e$ are positive integers satisfying $p\leqslant q$ and $pq-p<e<pq$, then for $G\in\mathcal{K}(p,q,e)$, we have
$$
\rho(G)\leqslant\rho(K_{p,q}^e),
$$
where $K_{p,q}^e\in \mathcal{K}(p,q,e)$ is the graph obtained from $K_{p,q}$ by deleting $pq-e$ edges incident with a common vertex in the partite set of order $q$.
\end{lemma}

\begin{remark}{\rm
In fact, from the proof of Lemma \ref{lem2.6}, we can deduce that the above equality holds if and only if $G\cong K_{p,q}^e$.}
\end{remark}

\begin{lemma}[\cite{BFP2008}]\label{lem2.7}
If $G$ is a bipartite graph, then
$\rho(G)\leqslant\sqrt{e(G)}.$
\end{lemma}

The following famous theorem of Folkman and Fulkerson \cite{FF1970} is one of the most important tools for our proof.
\begin{theorem}[\cite{FF1970}]\label{thm2.8}
Let $G=(X,Y)$ be a bipartite graph, and $g,f: V(G)\rightarrow \mathbb{Z}$ be functions such that $g(v)\leqslant f(v)$ for all $v\in V(G)$. Then $G$ has a $(g,f)$-factor if and only if
$$
\sum_{v\in S}f(v)+\sum_{v\in T}(d_G(v)-g(v))-e(S,T)\geqslant0
$$
and
$$
\sum_{v\in T}f(v)+\sum_{v\in S}(d_G(v)-g(v))-e(T,S)\geqslant0
$$
for all subsets $S\subseteq X$ and $T\subseteq Y.$
\end{theorem}
For two positive integers $a\leqslant b$, let $g(v)=a$ and $f(v)=b$ for all $v\in V(G)$ in Theorem \ref{thm2.8}, we obtain the following corollary, which gives a criterion for a bipartite graph to have an $[a,b]$-factor.
\begin{corollary}\label{cor2.9}
Let $G=(X,Y)$ be a bipartite graph, and $a\leqslant b$ be two positive integers. Then $G$ has an $[a,b]$-factor if and only if
$$
b|S|+\sum_{v\in T}d_G(v)-a|T|-e(S,T)\geqslant0
$$
and
$$
b|T|+\sum_{v\in S}d_G(v)-a|S|-e(T,S)\geqslant0
$$
for all subsets $S\subseteq X$ and $T\subseteq Y.$
\end{corollary}

The following lemma presents a property of bipartite graphs with $[a,b]$-factors.
\begin{lemma}\label{lem2.10}
Let $G=(X,Y)$ be a bipartite graph with $|X|\leqslant|Y|$, and $a\leqslant b$ be two positive integers. If $G$ has an $[a,b]$-factor, then $a|Y|\leqslant b|X|.$
\end{lemma}
\begin{proof}
Assume that $F$ is an $[a,b]$-factor of $G$, then $a\leqslant d_F(v)\leqslant b$ for all $v\in V(F)$. Furthermore, we have
$$
a|Y|\leqslant\sum_{v\in Y}d_F(v)=e(F)=\sum_{v\in X}d_F(v)\leqslant b|X|,
$$
as desired.
\end{proof}
\section{\normalsize Proofs of Theorems~\ref{thm1.1} and \ref{thm1.2}}\setcounter{equation}{0}
In this section, we give the proofs of Theorems~\ref{thm1.1} and \ref{thm1.2}. Before doing so, we need the following lemmas.

\begin{lemma}[\cite{WZ2023}]\label{lem3.1}
Let $a\leqslant b$ be two positive integers, and $G$ be a graph of order $n$ and $\delta(G)\geqslant a$. If
$$
e(G)\geqslant \dbinom{n-1}{2}+\frac{a+1}{2}
$$
and $na\equiv0\pmod{2}$ when $a=b$, then $G$ has an $[a,b]$-factor.
\end{lemma}

\begin{lemma}[\cite{OO1961}]\label{lem3.2}
Let $G$ be a graph of order $n\geqslant3$. If $e(G)\geqslant$
\begin{math}
\left( \begin{smallmatrix}
n-1\\2
\end{smallmatrix} \right),
\end{math}
then $G$ has a Hamilton path, unless $G\cong K_{n-1}\cup K_1$ or $G\cong K_{1,3}$.
\end{lemma}

\begin{lemma}[\cite{OO1961}]\label{lem3.3}
Let $G$ be a graph of order $n\geqslant3$. If $e(G)\geqslant$
\begin{math}
\left( \begin{smallmatrix}
n-1\\2
\end{smallmatrix} \right)+1,
\end{math}
then $G$ has a Hamilton cycle, unless $G\cong K_1\vee (K_{n-2}\cup K_1)$ or $G\cong K_2\vee \overline{K_3}$.
\end{lemma}

\begin{proof}[\bf Proof of Theorem \ref{thm1.1}]
Let $G$ be an $n$-vertex graph forbidding $[a,b]$-factors. We consider the following two possible cases.

{\bf Case 1.} $\delta(G)\leqslant a-1$.
In this case, there exists a vertex $u\in V(G)$ such that $d_G(u)\leqslant a-1$, which implies that $G$ is a spanning subgraph of $K_{a-1}\vee (K_{n-a}\cup K_1)$. Therefore, we have
$$
e(G)\leqslant e(K_{a-1}\vee (K_{n-a}\cup K_1))=\dbinom{n-1}{2}+a-1,
$$
where the first equality holds if and only if $G\cong K_{a-1}\vee (K_{n-a}\cup K_1)$.

{\bf Case 2.} $\delta(G)\geqslant a$.
In this case, since $G$ contains no $[a,b]$-factors, by Lemma \ref{lem3.1}, we have $e(G)<$
\begin{math}
\left( \begin{smallmatrix}
n-1\\2
\end{smallmatrix} \right)+\frac{a+1}{2},
\end{math}
i.e.,
$$
e(G)\leqslant \dbinom{n-1}{2}+\frac{a}{2}.
$$
When $a=1$, we have $e(G)\leqslant$
\begin{math}
\left( \begin{smallmatrix}
n-1\\2
\end{smallmatrix} \right)=\left( \begin{smallmatrix}
n-1\\2
\end{smallmatrix} \right)+a-1.
\end{math}
If
$e(G)=$
\begin{math}
\left( \begin{smallmatrix}
n-1\\2
\end{smallmatrix} \right),
\end{math}
then by Lemma \ref{lem3.2} and $\delta(G)\geqslant a=1$, we have $G\cong K_{1,3}$. Otherwise, $G$ has a Hamilton path. Combining with $na\equiv0\pmod{2}$ for $a=b$, it implies that $G$ has a $[1,b]$-factor, a contradiction. Note that if $b=1$ or $b=2$. Then $K_{1,3}$ contains no $[1,b]$-factors. If $b\geqslant3$, then $K_{1,3}$ is a $[1,b]$-factor of itself, a contradiction.
When $a=2$, we have $e(G)\leqslant$
\begin{math}
\left( \begin{smallmatrix}
n-1\\2
\end{smallmatrix} \right)+1=\left( \begin{smallmatrix}
n-1\\2
\end{smallmatrix} \right)+a-1.
\end{math}
Furthermore, if $e(G)=$
\begin{math}
\left( \begin{smallmatrix}
n-1\\2
\end{smallmatrix} \right)+1,
\end{math}
then by Lemma \ref{lem3.3} and $\delta(G)\geqslant a=2$, we have $G\cong K_2\vee \overline{K_3}$. Otherwise, $G$ has a Hamilton cycle, which implies that $G$ has a $[2,b]$-factor, a contradiction. Note that if $b=2$. Then $K_2\vee \overline{K_3}$ contains no $[2,b]$-factors. If $b\geqslant3$, then $K_2\vee \overline{K_3}$ is a $[2,b]$-factor of itself, a contradiction.
When $a\geqslant3$, we have
$e(G)\leqslant$
\begin{math}
\left( \begin{smallmatrix}
n-1\\2
\end{smallmatrix} \right)+\frac{a}{2}
\end{math}
$<$
\begin{math}
\left( \begin{smallmatrix}
n-1\\2
\end{smallmatrix} \right)+a-1.
\end{math}

This completes the proof.
\end{proof}

Next, we give the proof of Theorem \ref{thm1.2}.

\begin{proof}[\bf Proof of Theorem \ref{thm1.2}]
Suppose that $\rho(G)\geqslant \rho(K_{a-1}\vee (K_{n-a}\cup K_1))$ and $G\ncong K_{a-1}\vee (K_{n-a}\cup K_1)$. In order to obtain a contradiction, we will show that $G$ contains an $[a,b]$-factor.

We first assert that $G$ is connected. Note that $K_{a-1}\vee (K_{n-a}\cup K_1)$ contains $K_{n-1}\cup K_1$ as a spanning subgraph. By Corollary \ref{cor2.3}, we have $\rho(K_{a-1}\vee (K_{n-a}\cup K_1))\geqslant \rho(K_{n-1}\cup K_1)=n-2$, where the first equality holds if and only if $a=1$. If $G$ is not connected, then assume that $G_1,\ldots,G_h$ are the components of $G$. It is clear that $\rho(G)=\max\{\rho(G_1),\ldots,\rho(G_h)\}\leqslant\rho(K_{n-1}\cup K_1)\leqslant\rho(K_{a-1}\vee (K_{n-a}\cup K_1))$, where $\rho(G)=\rho(K_{a-1}\vee (K_{n-a}\cup K_1))$ holds if and only if $a=1$ and $G\cong K_{n-1}\cup K_1$, a contradiction to $G\ncong K_{a-1}\vee (K_{n-a}\cup K_1)$. Thus, $G$ is connected. Furthermore, by Lemma \ref{lem2.4} and $\rho(G)\geqslant \rho(K_{a-1}\vee (K_{n-a}\cup K_1))\geqslant n-2$, we have
$$
n-2\leqslant\rho(G)\leqslant\sqrt{2e(G)-n+1},
$$
which deduces that $e(G)\geqslant$
\begin{math}
\left( \begin{smallmatrix}
n-1\\2
\end{smallmatrix} \right)+1.
\end{math}
In what follows, assume that $e(G)=$
\begin{math}
\left( \begin{smallmatrix}
n-1\\2
\end{smallmatrix} \right)+t,
\end{math}
where $1\leqslant t\leqslant n-1$. We claim that $t\geqslant a$. If $1\leqslant t\leqslant a-1$, then by Lemma \ref{lem2.5}, we have $\rho(G)\leqslant \rho(K_{t}\vee (K_{n-1-t}\cup K_1))\leqslant\rho(K_{a-1}\vee (K_{n-a}\cup K_1))$, where $\rho(G)=\rho(K_{a-1}\vee (K_{n-a}\cup K_1))$ holds if and only if $G\cong K_{a-1}\vee (K_{n-a}\cup K_1)$, a contradiction to our assumption. Hence, $t\geqslant a$, which implies that $e(G)\geqslant$
\begin{math}
\left( \begin{smallmatrix}
n-1\\2
\end{smallmatrix} \right)+a.
\end{math}
By Theorem~\ref{thm1.1}, $G$ has an $[a,b]$-factor, as desired.

This completes the proof.
\end{proof}

\section{\normalsize Proofs of Theorems~\ref{thm1.3} and \ref{thm1.4}}\setcounter{equation}{0}
In this section, we give the proofs of Theorems~\ref{thm1.3} and \ref{thm1.4}. First, we prove Theorem~\ref{thm1.3}, which establishes an upper bound on the size of a bipartite graph with given partite sets forbidding $[a,b]$-factors.
\begin{proof}[\bf Proof of Theorem \ref{thm1.3}]
Let $G=(X,Y)$ be a bipartite graph forbidding $[a,b]$-factors. Assume that $|X|=p$ and $|Y|=q$, where $p\leqslant q$. Let $a,b$ be two positive integers with $a\leqslant b$.

(i)\ If $aq>bp$, then by Lemma \ref{lem2.10}, any bipartite graph with bipartite orders $p$ and $q$ contains no $[a,b]$-factors. Hence, $K_{p,q}$ contains no $[a,b]$-factors. Furthermore, $e(G)\leqslant e(K_{p,q})=pq$ with equality if and only if $G\cong K_{p,q}$, as desired.

(ii)\ If $a>p$, then for each vertex $v\in Y,$ we have $d_G(v)\leqslant p<a$, which implies that any bipartite graph with bipartite orders $p$ and $q$ contains no $[a,b]$-factors. Hence, $K_{p,q}$ contains no $[a,b]$-factors. Furthermore, $e(G)\leqslant e(K_{p,q})=pq$ with equality if and only if $G\cong K_{p,q}$.

(iii)\ We proceed with the following two possible cases.

{\bf Case 1.} $\delta(G)\leqslant a-1$.
In this case, there exists a vertex $u\in V(G)$ such that $d_G(u)\leqslant a-1$. If $u\in X$, then $G$ is a spanning subgraph of $D(p-1,1;a-1,q-a+1).$ Therefore,
\begin{align}
e(G)&\leqslant e(D(p-1,1;a-1,q-a+1))\label{eq:4.01}\\
&=(p-1)q+a-1\notag\\
&\leqslant p(q-1)+a-1.\label{eq:4.02}
\end{align}
Note that equality in \eqref{eq:4.01} holds if and only if $G\cong D(p-1,1;a-1,q-a+1)$; inequality in \eqref{eq:4.02} holds since $p\leqslant q$, and equality in \eqref{eq:4.02} holds if and only if $p=q$. Hence, $e(G)=p(q-1)+a-1$ holds if and only if $G\cong D(p-1,1;a-1,p-a+1)$, as desired.
If $u\in Y,$ then $G$ is a spanning subgraph of $D(a-1,p-a+1;q-1,1).$ Therefore,
$$
e(G)\leqslant e(D(a-1,p-a+1;q-1,1))=p(q-1)+a-1,
$$
where the first equality holds if and only if $G\cong D(a-1,p-a+1;q-1,1)$.

{\bf Case 2.} $\delta(G)\geqslant a$.
In this case, if $b\geqslant q$, then for each vertex $v\in V(G)$, we have $a\leqslant d_G(v)\leqslant q\leqslant b$, which implies that $G$ is an $[a,b]$-factor of itself, a contradiction. If $b<q$, then in order to complete the proof, it suffices to prove that $e(G)<p(q-1)+a-1$.
Since $G$ contains no $[a,b]$-factors, by Corollary \ref{cor2.9}, there exist two vertex subsets $S\subseteq X$ and $T\subseteq Y$ such that
$$
\gamma^*(S,T):=b|S|+\sum_{v\in T}d_G(v)-a|T|-e(S,T)<0
$$
or
$$
\gamma^*(T,S):=b|T|+\sum_{v\in S}d_G(v)-a|S|-e(T,S)<0.
$$
Choose such a pair $(S,T)$ so that $S\cup T$ is maximal. Then we proceed by distinguishing the following two possible subcases.

{\bf Subcase 2.1.} $\gamma^*(S,T)<0$.
In this subcase, we have
\begin{align}\label{eq:3.03}
e(X\backslash S,T)=\sum_{v\in T}d_{G-S}(v)=\sum_{v\in T}d_G(v)-e(S,T)<a|T|-b|S|.
\end{align}
Moreover, we have the following claims.
\begin{claim}\label{3}
$|T|\geqslant b+1$.
\end{claim}
\begin{proof}[\bf Proof of Claim~\ref{3}]
Suppose that $|T|\leqslant b$. Together with $\delta(G)\geqslant a$, we have
$$
\gamma^*(S,T)=b|S|+\sum_{v\in T}d_G(v)-a|T|-e(S,T)\geqslant b|S|-e(S,T)\geqslant b|S|-|S||T|=(b-|T|)|S|\geqslant0,
$$
a contradiction. So Claim~\ref{3} holds.
\end{proof}

\begin{claim}\label{5}
If $T\subset Y,$ then for each $v\in Y\backslash T$, $e(X\backslash S,\{v\})>a$.
\end{claim}
\begin{proof}[\bf Proof of Claim~\ref{5}]
Suppose that there exists a vertex $u\in Y\backslash T$ such that $e(X\backslash S,\{u\})\leqslant a$. Let $T'=T\cup \{u\}$. Then
\begin{align*}
\gamma^*(S,T')&=b|S|+\sum_{v\in T'}d_G(v)-a|T'|-e(S,T')\\
&=b|S|+\sum_{v\in T}d_G(v)+d_G(u)-a(|T|+1)-(e(S,T)+e(S,\{u\}))\\
&=\gamma^*(S,T)+d_G(u)-a-e(S,\{u\})\\
&=\gamma^*(S,T)+e(X\backslash S,\{u\})-a<0,
\end{align*}
which contradicts the choice of $(S,T)$. Hence, Claim~\ref{5} holds.
\end{proof}

Recall that $|X|=p$ and $|Y|=q$. It follows from Claim~\ref{5} that if $T\subset Y,$ then for each $v\in Y\backslash T$,
\begin{align}\label{eq:3.04}
a<e(X\backslash S,\{v\})\leqslant|X\backslash S|=|X|-|S|=p-|S|.
\end{align}

\begin{claim}\label{6}
If $p=q$, then $T\subset Y.$
\end{claim}
\begin{proof}[\bf Proof of Claim~\ref{6}]
Suppose that $T=Y.$
Applying \eqref{eq:3.03} for $T=Y,$ we have $e(X\backslash S,Y)<a|Y|-b|S|$. Note that $e(X\backslash S,Y)=\sum_{v\in X\backslash S}d_{G}(v)$. Combining with $\delta(G)\geqslant a$, we obtain
$$a(p-|S|)=a(|X|-|S|)\leqslant\sum_{v\in X\backslash S}d_{G}(v)=e(X\backslash S,Y)<a|Y|-b|S|\leqslant a(|Y|-|S|)=a(q-|S|),$$
a contradiction. So Claim~\ref{6} holds.
\end{proof}

In what follows, we consider the size of $G$ according to whether $T\subset Y$ or not.
If $T\subset Y,$ then we have
\begin{align}
e(G)&=e(S,T)+e(X\backslash S,T)+e(X,Y\backslash T)\notag\\
&\leqslant|S||T|+e(X\backslash S,T)+p(q-|T|)\notag\\
&<|S||T|+a|T|-b|S|+p(q-|T|)\tag{by \eqref{eq:3.03}}\\
&=(|S|+a-p)|T|-b|S|+pq\notag\\
&\leqslant(|S|+a-p)(b+1)-b|S|+pq\tag{by \eqref{eq:3.04} and Claim~\ref{3}}\\
&=p(q-1)-b(p-a)+|S|+a\notag\\
&\leqslant p(q-1)-b(p-a)+p-a-1+a\tag{by \eqref{eq:3.04}}\\
&=p(q-1)-(b-1)(p-a)+a-1\notag\\
&\leqslant p(q-1)+a-1,\tag{by $b\geqslant1$ and \eqref{eq:3.04}}\notag
\end{align}
as desired.
If $T=Y,$ then we first assert that $S\subset X$. Suppose not, then $S=X$, which gives $\gamma^*(S,T)=\gamma^*(X,Y)$. Note that $\gamma^*(S,T)<0$. Therefore,
$$
\gamma^*(X,Y)=b|X|+\sum_{v\in Y}d_G(v)-a|Y|-e(X,Y)=b|X|-a|Y|=bp-aq<0.
$$
This implies that $bp<aq$, a contradiction.
Hence, $S\subset X$ and so $|S|\leqslant |X|-1=p-1$. Furthermore, we have
\begin{align}
e(G)&=e(S,Y)+e(X\backslash S,Y)\notag\\
&\leqslant|S||Y|+e(X\backslash S,Y)\notag\\
&<|S||Y|+a|Y|-b|S|\tag{by \eqref{eq:3.03}}\\
&=|S|(|Y|-b)+a|Y|\notag\\
&\leqslant(p-1)(q-b)+aq\tag{since $|S|\leqslant p-1$}\\
&=pq-bp-q+b+aq\notag\\
&=p(q-1)-(bp-aq)-(q-p)+b.\label{eq:3.05}
\end{align}
Note that $bp\geqslant aq$ and $q>p$ for $T=Y$ by Claim~\ref{6}.
If $bp-aq\geqslant b-a$, then in view of \eqref{eq:3.05}, we have
$$
e(G)<p(q-1)-(bp-aq)-(q-p)+b\leqslant p(q-1)-(b-a)-1+b=p(q-1)+a-1,
$$
as desired. If $0\leqslant bp-aq\leqslant b-a-1$, then $p\leqslant \frac{b-a-1+aq}{b}$, which implies that
$$
q-p\geqslant q-\frac{b-a-1+aq}{b}=\frac{(b-a)(q-1)+1}{b}>b-a,
$$
where the last inequality follows from $b<q$. Furthermore, in view of \eqref{eq:3.05}, we obtain
$$
e(G)<p(q-1)-(bp-aq)-(q-p)+b\leqslant p(q-1)-(b-a+1)+b=p(q-1)+a-1,
$$
as required.

{\bf Subcase 2.2.} $\gamma^*(T,S)<0$. In this subcase, we have
\begin{align}\label{eq:3.06}
e(S,Y\backslash T)=\sum_{v\in S}d_{G-T}(v)=\sum_{v\in S}d_G(v)-e(T,S)<a|S|-b|T|.
\end{align}
Moreover, we have the following claims.
\begin{claim}\label{7}
$|S|\geqslant b+1$.
\end{claim}
\begin{proof}[\bf Proof of Claim~\ref{7}]
Suppose that $|S|\leqslant b$. Then
\begin{align}\label{eq:3.008}
\gamma^*(T,S)&=b|T|+\sum_{v\in S}d_G(v)-a|S|-e(T,S)\notag\\
             &\geqslant b|T|-e(T,S)\\
             &\geqslant b|T|-|T||S|\notag\\
             &=(b-|S|)|T|\notag\\
             &\geqslant0,\notag
\end{align}
where the inequality in \eqref{eq:3.008} follows from $\delta(G)\geqslant a$, a contradiction. Hence, Claim~\ref{7} holds.
\end{proof}

\begin{claim}\label{8}
$S\subset X$.
\end{claim}
\begin{proof}[\bf Proof of Claim~\ref{8}]
If not, then $S=X$. Applying \eqref{eq:3.06} for $S=X$ and by $\delta(G)\geqslant a$, we have
$$
a(q-|T|)=a(|Y|-|T|)\leqslant\sum_{v\in Y\backslash T}d_{G}(v)=e(X,Y\backslash T)<a|X|-b|T|\leqslant a(|X|-|T|)=a(p-|T|),
$$
a contradiction. Hence, Claim~\ref{8} holds.
\end{proof}

It follows immediately from Claim~\ref{8} that $X\backslash S\neq\emptyset$. For each $v\in X\backslash S$, we have the following property.

\begin{claim}\label{9}
For each $v\in X\backslash S$, we have $e(Y\backslash T,\{v\})>a$.
\end{claim}
\begin{proof}[\bf Proof of Claim~\ref{9}]
Suppose that there exists a vertex $u\in X\backslash S$ such that $e(Y\backslash T,\{u\})\leqslant a$. Let $S'=S\cup \{u\}$. Then
\begin{align*}
\gamma^*(T,S')&=b|T|+\sum_{v\in S'}d_G(v)-a|S'|-e(T,S')\\
&=b|T|+\sum_{v\in S}d_G(v)+d_G(u)-a(|S|+1)-(e(T,S)+e(T,\{u\}))\\
&=\gamma^*(T,S)+d_G(u)-a-e(T,\{u\})\\
&=\gamma^*(T,S)+e(Y\backslash T,\{u\})-a<0,
\end{align*}
which contradicts the choice of $(S,T)$. So Claim~\ref{9} holds.
\end{proof}
By Claim~\ref{9}, for each vertex $v\in X\backslash S$, we have
\begin{align}\label{eq:3.07}
a<e(Y\backslash T,\{v\})\leqslant|Y\backslash T|=|Y|-|T|=q-|T|.
\end{align}

In what follows, we consider the size of $G$. By a direct computation, we obtain the size of $G$ as
\begin{align}
e(G)&=e(S,T)+e(S,Y\backslash T)+e(X\backslash S,Y)\notag\\
&\leqslant|S||T|+e(S,Y\backslash T)+(p-|S|)q\notag\\
&<|S||T|+a|S|-b|T|+(p-|S|)q\tag{by \eqref{eq:3.06}}\\
&=|S|(|T|+a-q)-b|T|+pq\notag\\
&\leqslant(b+1)(|T|+a-q)-b|T|+pq\tag{by Claim~\ref{7} and \eqref{eq:3.07}}\\
&=pq-q-b(q-a)+|T|+a\notag\\
&\leqslant p(q-1)-b(q-a)+|T|+a\tag{since $p\leqslant q$}\\
&\leqslant p(q-1)-b(q-a)+q-a-1+a\tag{by \eqref{eq:3.07}}\\
&=p(q-1)-(b-1)(q-a)+a-1\notag\\
&\leqslant p(q-1)+a-1,\tag{by $b\geqslant1$ and \eqref{eq:3.07}}
\end{align}
as desired.

Combining with Case 1 and Case 2, we complete the proof.
\end{proof}

Based on Theorem \ref{thm1.3}, we give the proof of Theorem~\ref{thm1.4}, which establishes an upper bound on the size of an $n$-vertex bipartite graph forbidding $[a,b]$-factors.

\begin{proof}[\bf Proof of Theorem~\ref{thm1.4}]
Let $G=(X,Y)$ be a bipartite graph of order $n$ forbidding $[a,b]$-factors. Assume that $|X|=p$ and $|Y|=q$, where $p\leqslant q$ and $p+q=n$. Let $a,b$ be two positive integers such that $a\leqslant b$ and $a\leqslant\lfloor\frac{n}{2}\rfloor$.

To prove our theorem, it suffices to show:
\begin{wst}
\item[{\rm (1)}]If $aq>bp$, then $e(G)\leqslant \lfloor\frac{an-1}{a+b}\rfloor(n-\lfloor\frac{an-1}{a+b}\rfloor)$ with equality if and only if $G\cong K_{\lfloor\frac{an-1}{a+b}\rfloor,n-\lfloor\frac{an-1}{a+b}\rfloor}$.
\item[{\rm (2)}]If $aq\leqslant bp$, then $e(G)\leqslant
    \lfloor\frac{n}{2}\rfloor(\lceil\frac{n}{2}\rceil-1)+a-1$ with equality if and only if one of the following holds:
\begin{wst}
\item[{\rm (i)}]$G\cong D(a-1,\frac{n}{2}-a;\frac{n}{2},1)$ or $D(a-1,\frac{n}{2}-a+1;\frac{n}{2}-1,1)$ for even $n$;
\item[{\rm (ii)}]$G\cong D(a-1,\frac{n+1}{2}-a;\frac{n-1}{2},1)$ for odd $n$.
\end{wst}
\end{wst}

We first give the proof of (1). If $aq>bp$, then by Theorem \ref{thm1.3}(i), we have $e(G)\leqslant pq$ with equality if and only if $G\cong K_{p,q}$. In addition, if $aq>bp$, then $aq-1\geqslant bp$, i.e., $a(n-p)-1\geqslant bp$, which implies that $p\leqslant\lfloor\frac{an-1}{a+b}\rfloor$. Note that $pq=p(n-p)\leqslant \lfloor\frac{an-1}{a+b}\rfloor(n-\lfloor\frac{an-1}{a+b}\rfloor)$, where equality holds if and only if $p=\lfloor\frac{an-1}{a+b}\rfloor$. Hence, $e(G)\leqslant \lfloor\frac{an-1}{a+b}\rfloor(n-\lfloor\frac{an-1}{a+b}\rfloor)$ with equality if and only if $G\cong K_{\lfloor\frac{an-1}{a+b}\rfloor,n-\lfloor\frac{an-1}{a+b}\rfloor}$.

Next we prove (2). We proceed by distinguishing the following two possible cases.

{\bf Case 1.} $a\leqslant p$.
In this case, by Theorem~\ref{thm1.3}(iii), we have $e(G)\leqslant p(n-p-1)+a-1$ with equality if and only if $G\cong D(a-1,p-a+1;n-p-1,1)$. Note that
$$
p(n-p-1)+a-1\leqslant \left\lfloor\frac{n}{2}\right\rfloor(\left\lceil\frac{n}{2}\right\rceil-1)+a-1,
$$
with equality if and only if $p=\frac{n}{2}-1$ or $p=\frac{n}{2}$ for even $n$, or $p=\frac{n-1}{2}$ for odd $n$. Hence, $e(G)\leqslant \lfloor\frac{n}{2}\rfloor(\lceil\frac{n}{2}\rceil-1)+a-1$ with equality if and only if $G\cong D(a-1,\frac{n}{2}-a;\frac{n}{2},1)$ or $G\cong D(a-1,\frac{n}{2}-a+1;\frac{n}{2}-1,1)$ for even $n$, or $G\cong D(a-1,\frac{n+1}{2}-a;\frac{n-1}{2},1)$ for odd $n$, as desired.

{\bf Case 2.} $a>p$.
In this case, by Theorem~\ref{thm1.3}(ii), one sees
\begin{align}
e(G)&\leqslant e(K_{p,n-p})=p(n-p)\notag\\
&\leqslant (a-1)(n-a+1)\label{eq:4.09}\\
&=(a-1)(n-a)+a-1\notag\\
&\leqslant (\left\lfloor\frac{n}{2}\right\rfloor-1)(n-\left\lfloor\frac{n}{2}\right\rfloor)+a-1\label{eq:4.10}\\
&=(\left\lfloor\frac{n}{2}\right\rfloor-1)\left\lceil\frac{n}{2}\right\rceil+a-1\notag\\
&\leqslant\left\lfloor\frac{n}{2}\right\rfloor(\left\lceil\frac{n}{2}\right\rceil-1)+a-1\label{eq:4.11}.
\end{align}
Note that \eqref{eq:4.09} follows from $p\leqslant a-1$, and equality in \eqref{eq:4.09} holds if and only if $p=a-1$; \eqref{eq:4.10} follows from $a\leqslant \lfloor\frac{n}{2}\rfloor$, and equality in \eqref{eq:4.10} holds if and only if $a=\lfloor\frac{n}{2}\rfloor$; \eqref{eq:4.11} follows from $\lfloor\frac{n}{2}\rfloor\leqslant \lceil\frac{n}{2}\rceil$, and equality in \eqref{eq:4.11} holds if and only if $\lfloor\frac{n}{2}\rfloor=\lceil\frac{n}{2}\rceil$, i.e., $n$ is even. Hence, $e(G)=\lfloor\frac{n}{2}\rfloor(\lceil\frac{n}{2}\rceil-1)+a-1$ holds if and only if $G\cong K_{p,n-p}$, $p=a-1$, $a=\lfloor\frac{n}{2}\rfloor$ and $n$ is even, which is equivalent to $n$ is even, $a=\frac{n}{2}$ and $G\cong K_{\frac{n}{2}-1,\frac{n}{2}+1}$, as desired.
\end{proof}

\section{\normalsize Proofs of Theorems~\ref{thm1.5} and \ref{thm1.7}}\setcounter{equation}{0}
In this section, we shall give the proofs of Theorems~\ref{thm1.5} and \ref{thm1.7}. Firstly, we prove Theorem~\ref{thm1.5}, which provides an upper bound on the spectral radius of a bipartite graph $G$ with given partite sets forbidding $[a,b]$-factors.
\begin{proof}[\bf Proof of Theorem \ref{thm1.5}]
Let $G=(X,Y)$ be a bipartite graph forbidding $[a,b]$-factors. Assume that $|X|=p$ and $|Y|=q$, where $p\leqslant q$. Let $a,b$ be two positive integers with $a\leqslant b$.

(i)\ \ If $aq>bp$, then by Lemma \ref{lem2.10}, any bipartite graph with bipartite orders $p$ and $q$ contains no $[a,b]$-factors. Hence, $K_{p,q}$ contains no $[a,b]$-factors. Furthermore, $\rho(G)\leqslant \rho(K_{p,q})=\sqrt{pq}$ with equality if and only if $G\cong K_{p,q}$.

(ii)\ \ If $a>p$, then for each vertex $v\in Y$, we have $d_G(v)\leqslant p<a$, which implies that any bipartite graph with bipartite orders $p$ and $q$ contains no $[a,b]$-factors. Hence, $K_{p,q}$ contains no $[a,b]$-factors. Furthermore, $\rho(G)\leqslant \rho(K_{p,q})=\sqrt{pq}$ with equality if and only if $G\cong K_{p,q}$.

(iii)\ \ Suppose that $\rho(G)\geqslant\rho(D(a-1,p-a+1;q-1,1))$ and $G\ncong D(a-1,p-a+1;q-1,1)$. In what follows, we are to show that $G$ contains an $[a,b]$-factor.
For convenience, let $G_1:=D(a-1,p-a+1;q-1,1)$. Then we have the following claim.
\begin{claim}\label{2}
$e(G)\geqslant p(q-1)$ if $a=1$, and $e(G)>p(q-1)$ if $a>1$.
\end{claim}
\begin{proof}[\bf Proof of Claim~\ref{2}]
Note that $a$ is a positive integer. If $a=1$, then $G_1=K_{p,q-1}\cup K_1$. Thus, $\rho(G)\geqslant\rho(K_{p,q-1}\cup K_1)=\sqrt{p(q-1)}$. Together with Lemma \ref{lem2.7}, we have
$$
\sqrt{p(q-1)}\leqslant\rho(G)\leqslant\sqrt{e(G)},
$$
which implies that $e(G)\geqslant p(q-1)$.
If $a>1$, then $G_1$ is connected and contains $K_{p,q-1}\cup K_1$ as a proper spanning subgraph. By Lemma \ref{lem2.2}, we have $\rho(G_1)>\rho(K_{p,q-1}\cup K_1)=\sqrt{p(q-1)}$. Combining with $\rho(G)\geqslant\rho(G_1)$ and Lemma \ref{lem2.7}, we obtain
$$
\sqrt{p(q-1)}<\rho(G_1)\leqslant\rho(G)\leqslant\sqrt{e(G)},
$$
which deduces that $e(G)>p(q-1)$.
Hence, Claim~\ref{2} holds.
\end{proof}

By Claim~\ref{2}, we may assume that $e(G)=p(q-1)+r$, where $0\leqslant r\leqslant p$.
Next we assert that $r=0$ or $r\geqslant a$. If not, then $1\leqslant r\leqslant a-1$. By Lemma \ref{lem2.6}, we have
$$
\rho(G)\leqslant \rho(D(r,p-r;q-1,1))
$$
with equality if and only if $G\cong D(r,p-r;q-1,1)$. It is easy to see that $D(r,p-r;q-1,1)$ is a spanning subgraph of $D(a-1,p-a+1;q-1,1)$. Thus, by Lemma \ref{lem2.2}, $\rho(D(r,p-r;q-1,1))\leqslant\rho(D(a-1,p-a+1;q-1,1))$ with equality if and only if $r=a-1$. Hence, $\rho(G)\leqslant\rho(D(a-1,p-a+1;q-1,1))$ with equality if and only if $G\cong D(a-1,p-a+1;q-1,1)$, which contradicts our assumption. Therefore, we have $r=0$ or $r\geqslant a$. If $r=0$, then $e(G)=p(q-1)$. By Claim~\ref{2}, we have $a=1$. Furthermore, by Theorem \ref{thm1.3}(iii), $G$ contains a $[1,b]$-factor, as required. If $r\geqslant a$, then $e(G)\geqslant p(q-1)+a$. By Theorem \ref{thm1.3}(iii), $G$ contains an $[a,b]$-factor, as desired.

This completes the proof.
\end{proof}

Based on Theorem~\ref{thm1.5}, we give the proof of Theorem~\ref{thm1.7}, which establishes an upper bound on the spectral radius of a bipartite graph with given order forbidding $[a,b]$-factors.
\begin{proof}[\bf Proof of Theorem~\ref{thm1.7}]
Let $G=(X,Y)$ be a bipartite graph of order $n$ forbidding $[a,b]$-factors and assume that $|X|=p$ and $|Y|=q$, where $p\leqslant q$ and $p+q=n$. Let $a,b$ be two positive integers with $a\leqslant b$ and $a\leqslant \lfloor\frac{n}{2}\rfloor$.

To prove our theorem, it suffices to show:
\begin{wst}
\item[{\rm (1)}]If $aq>bp$, then $\rho(G)\leqslant \sqrt{\lfloor\frac{an-1}{a+b}\rfloor(n-\lfloor\frac{an-1}{a+b}\rfloor)}$ with equality if and only if $G\cong K_{\lfloor\frac{an-1}{a+b}\rfloor,n-\lfloor\frac{an-1}{a+b}\rfloor}$.
\item[{\rm (2)}]If $aq\leqslant bp$, then
$$
\rho(G)\leqslant \rho(D(a-1,\left\lceil\frac{n}{2}\right\rceil-a;\left\lfloor\frac{n}{2}\right\rfloor,1))
$$
with equality if and only if $G\cong D(a-1,\left\lceil\frac{n}{2}\right\rceil-a;\left\lfloor\frac{n}{2}\right\rfloor,1)$.
\end{wst}

We first prove (1). If $aq>bp$, then by Theorem \ref{thm1.5}(i), we have $\rho(G)\leqslant \sqrt{pq}$ with equality if and only if $G\cong K_{p,q}$. In addition, if $aq>bp$, then $aq-1\geqslant bp$, i.e., $a(n-p)-1\geqslant bp$, which implies that $p\leqslant\lfloor\frac{an-1}{a+b}\rfloor$. Note that $\sqrt{pq}=\sqrt{p(n-p)}\leqslant \sqrt{\lfloor\frac{an-1}{a+b}\rfloor(n-\lfloor\frac{an-1}{a+b}\rfloor)}$ with equality if and only if $p=\lfloor\frac{an-1}{a+b}\rfloor$. Hence, $\rho(G)\leqslant \sqrt{\lfloor\frac{an-1}{a+b}\rfloor(n-\lfloor\frac{an-1}{a+b}\rfloor)}$ with equality if and only if $G\cong K_{\lfloor\frac{an-1}{a+b}\rfloor,n-\lfloor\frac{an-1}{a+b}\rfloor}$.

Next we prove (2). We consider the following two possible cases.

{\bf Case 1.} $a\leqslant p$.
In this case, by Theorem~\ref{thm1.5}(iii), we have $\rho(G)\leqslant\rho(D(a-1,p-a+1;q-1,1))$ with equality if and only if $G\cong D(a-1,p-a+1;q-1,1)$. In what follows, we tend to show
$$
\rho(D(a-1,p-a+1;q-1,1))\leqslant \rho(D(a-1,\frac{n}{2}-a;\frac{n}{2},1))
$$
for even $n$, where equality holds if and only if $p=\frac{n}{2}-1$, or $p=\frac{n}{2}$ and $a=1$, and
$$
\rho(D(a-1,p-a+1;q-1,1))\leqslant \rho(D(a-1,\frac{n+1}{2}-a;\frac{n-1}{2},1))
$$
for odd $n$, where equality holds if and only if $p=\frac{n-1}{2}$.

We first consider that $n$ is even.
If $p=\frac{n}{2}-1$, then the result holds obviously. In what follows,
it suffices to show that $\rho(D(a-1,p-a+1;q-1,1))\leqslant\rho(D(a-1,\frac{n}{2}-a;\frac{n}{2},1))$ for $p\neq\frac{n}{2}-1$, where equality holds if and only if $p=\frac{n}{2}$ and $a=1$.

For convenience, let $G_1:=D(a-1,p-a+1;q-1,1)$ and $G_2:=D(a-1,\frac{n}{2}-a;\frac{n}{2},1)$. If $a=1$, then $G_1=K_{p,q-1}\cup K_1$ and $G_2=K_{\frac{n}{2}-1,\frac{n}{2}}\cup K_1$. It is easy to check that
$$
\sqrt{p(q-1)}=\rho(G_1)\leqslant\rho(G_2)=\sqrt{\frac{n}{2}(\frac{n}{2}-1)}
$$
with equality if and only if $p=\frac{n}{2}-1$ or $p=\frac{n}{2}$.

Next we assume that $a\geqslant2$.
Let $X'$ and $Y'$ be two partite sets of $G_1$, where $|X'|=p$ and $|Y'|=q$. Furthermore, let
$$
X'_1=\{v\in X':d_{G_1}(v)=q\},\ \ X'_2=\{v\in X':d_{G_1}(v)=q-1\},
$$
and let
$$
Y'_1=\{v\in Y':d_{G_1}(v)=p\},\ \ Y'_2=\{v\in Y':d_{G_1}(v)=a-1\}.
$$
Consider the partition $\pi_1:\ V(G_1)=X'_1\cup X'_2\cup Y'_1\cup Y'_2$. Then the corresponding quotient matrix of $A(G_1)$ is
$$
M_1=\left(
  \begin{array}{cccc}
    0 & 0 & q-1 & 1\\
    0 & 0 & q-1 & 0\\
    a-1 & p-a+1 & 0 & 0\\
    a-1 & 0 & 0 & 0\\
  \end{array}
\right).
$$
By a simple calculation, we obtain the characteristic polynomial of $M_1$ as
\begin{align}\label{eq:3.01}
\Phi_1(x)=x^4-(pq-p+a-1)x^2+(a-1)(q-1)(p-a+1).
\end{align}
Note that the partition $\pi_1$ is equitable. Hence, by Lemma \ref{lem2.1}, the largest root of $\Phi_1(x)=0$ equals the spectral radius of $G_1.$

Replacing $p$ with $\frac{n}{2}-1$ in \eqref{eq:3.01}, we have
\begin{align}\label{eq:3.02}
\Phi_2(x)=\frac{1}{4}(4x^4-(n^2-2n+4a-4)t^2+n(a-1)(n-2a)).
\end{align}
It is clear that the largest root of $\Phi_2(x)=0$ equals the spectral radius of $G_2.$

In view of \eqref{eq:3.01}, \eqref{eq:3.02} and by a direct computation, we have
$$
\Phi_2(x)-\Phi_1(x)=-\frac{1}{4}(n-2p-2)((n-2p)x^2-(a-1)(n-2p+2a-2)).
$$
Let $f_1(x)=(n-2p)x^2-(a-1)(n-2p+2a-2)$ be a real function in $x$. If $p=\frac{n}{2}$, then $f_1(x)=-2(a-1)^2<0$ for $a\geqslant2$. Thus,  $\Phi_2(x)-\Phi_1(x)=-\frac{1}{4}(n-2p-2)f_1(x)=-(a-1)^2<0$. It follows that $\Phi_2(\rho(G_1))-\Phi_1(\rho(G_1))<0$. Since $\Phi_1(\rho(G_1))=0$, we have $\Phi_2(\rho(G_1))<0$, which implies that
$\rho(G_1)<\rho(G_2)$, as desired.
If $p\leqslant \frac{n}{2}-2$, then $n-2p>0$. Consider the derivative of $f_1(x)$, we have $f'_1(x)=2(n-2p)x>0$ for $x>0$, which deduces that $f_1(x)$ is a monotonically increasing function for $x>0$. Note that $K_{a-1,q}$ is a proper subgraph of $G_1$. Then by Lemma \ref{lem2.2}, we have $\rho(G_1)>\rho(K_{a-1,q})=\sqrt{(a-1)q}$.
Hence,
$$
f_1(\rho(G_1))>f_1(\sqrt{(a-1)q})=(a-1)(2p^2-(3n-2)p+n^2-n-2a+2).
$$
Let $f_2(x)=2x^2-(3n-2)x+n^2-n-2a+2$ be a real function in $x$, where $x\leqslant \frac{n}{2}-2$. Then $f'_2(x)=4x-3n+2\leqslant4(\frac{n}{2}-2)-3n+2=-n-6<0$, which implies that $f_2(x)$ is a monotonically decreasing function for $x\leqslant \frac{n}{2}-2$. Thus, $f_2(x)\geqslant f_2(\frac{n}{2}-2)=2n-2a+6>0$. Hence, $f_2(p)>0$, which deduces that $f_1(\rho(G_1))>f_1(\sqrt{(a-1)q})=(a-1)f_2(p)>0$. Furthermore, we have $\Phi_2(\rho(G_1))-\Phi_1(\rho(G_1))=-\frac{1}{4}(n-2p-2)f_1(\rho(G_1))<0$. Since $\Phi_1(\rho(G_1))=0$, we have $\Phi_{2}(\rho(G_1))<0$, which implies that $\rho(G_1)<\rho(G_2)$, as desired.

Next we consider that $n$ is odd.
If $p=\frac{n-1}{2}$, then the result holds obviously. In what follows,
it suffices to show that $\rho(D(a-1,p-a+1;q-1,1))<\rho(D(a-1,\frac{n+1}{2}-a;\frac{n-1}{2},1))$ for $p<\frac{n-1}{2}$.

For convenience, let $G_3:=D(a-1,\frac{n+1}{2}-a;\frac{n-1}{2},1)$. If $a=1$, then $G_1=K_{p,q-1}\cup K_1$ and $G_3=K_{\frac{n-1}{2},\frac{n-1}{2}}\cup K_1$. It is easy to check that
$$
\sqrt{p(q-1)}=\rho(G_1)\leqslant\rho(G_3)=\frac{n-1}{2}
$$
with equality if and only if $p=\frac{n-1}{2}$.

Next we assume that $a\geqslant2$. Replacing $p$ with $\frac{n-1}{2}$ in \eqref{eq:3.01}, we have
\begin{align}\label{eq:5.04}
\Phi_3(x)=\frac{1}{4}(4x^4-(n^2-2n+4a-3)t^2+(a-1)(n-1)(n-2a+1)).
\end{align}
It is clear that the largest root of $\Phi_3(x)=0$ equals the spectral radius of $G_3.$

In view of \eqref{eq:3.01}, \eqref{eq:5.04} and by a direct calculation, we have
$$
\Phi_3(x)-\Phi_1(x)=-\frac{1}{4}(n-2p-1)((n-2p-1)x^2-(a-1)(n-2p+2a-3)).
$$
Let $f_3(x)=(n-2p-1)x^2-(a-1)(n-2p+2a-3)$ be a real function in $x$ with $x>0$. Then $f'_3(x)=2(n-2p-1)x>0$ for $p<\frac{n-1}{2}$, which implies that $f_3(x)$ is a monotonically increasing function for $x>0$. Note that $K_{a-1,q}$ is a proper subgraph of $G_1$. Then by Lemma \ref{lem2.2}, we have $\rho(G_1)>\rho(K_{a-1,q})=\sqrt{(a-1)q}$. Hence,
$$
f_3(\rho(G_1))>f_3(\sqrt{(a-1)q})=(a-1)(2p^2-(3n-3)p+n^2-2n-2a+3).
$$
Let $f_4(x)=2x^2-(3n-3)x+n^2-2n-2a+3$ be a real function in $x$, where $x\leqslant \frac{n-3}{2}$. Then $f'_4(x)=4x-3n+3\leqslant4\times\frac{n-3}{2}-3n+3=-n-3<0$, which implies that $f_4(x)$ is a monotonically decreasing function for $x\leqslant \frac{n-3}{2}$. Thus, $f_4(x)\geqslant f_4(\frac{n-3}{2})=n-2a+3>0$. Hence, $f_4(p)>0$, which deduces that $f_3(\rho(G_1))>f_3(\sqrt{(a-1)q})=(a-1)f_4(p)>0$. Furthermore, we have $\Phi_3(\rho(G_1))-\Phi_1(\rho(G_1))=-\frac{1}{4}(n-2p-1)f_3(\rho(G_1))<0$. Since $\Phi_1(\rho(G_1))=0$, we have $\Phi_{3}(\rho(G_1))<0$, which implies that $\rho(G_1)<\rho(G_3)$, as desired.

{\bf Case 2.} $a>p$.
In this case, by Theorem~\ref{thm1.5}(ii), we have $\rho(G)\leqslant \sqrt{pq}\leqslant\linebreak\sqrt{(a-1)(n-a+1)}$, where $\rho(G)=\sqrt{(a-1)(n-a+1)}$ holds if and only if $G\cong K_{a-1,n-a+1}$.
In order to complete the proof, it suffices to show that
$$
\sqrt{(a-1)(n-a+1)}\leqslant\rho(D(a-1,\frac{n}{2}-a;\frac{n}{2},1))=\rho(G_2)
$$
for even $n$, where equality holds if and only if $a=\frac{n}{2}$, and
$$
\sqrt{(a-1)(n-a+1)}<\rho(D(a-1,\frac{n+1}{2}-a;\frac{n-1}{2},1))=\rho(G_3)
$$
for odd $n$.

We first consider that $n$ is even. If $a=\frac{n}{2}$, then the result holds obviously. Thus, in what follows, it suffices to prove
$\sqrt{(a-1)(n-a+1)}<\rho(G_2)$ for $1\leqslant a<\frac{n}{2}$.

If $a=1$, then $G_2=K_{\frac{n}{2}-1,\frac{n}{2}}\cup K_1$. It is easy to check that
$$
0=\sqrt{(a-1)(n-a+1)}<\sqrt{\frac{n}{2}(\frac{n}{2}-1)}=\rho(G_2),
$$
as desired.

Next we consider $2\leqslant a<\frac{n}{2}$. Recall that $\rho(G_2)$ equals the largest root of $\Phi_2(x)=0$, where $\Phi_2(x)$ is given in \eqref{eq:3.02}.
Substituting $x=\sqrt{(a-1)(n-a+1)}$ into \eqref{eq:3.02} yields
$$
\Phi_2(\sqrt{(a-1)(n-a+1)})=-\frac{1}{4}(a-1)(n-2a)(2a^2-(3n+4)a+n^2+2n+2).
$$
Let $g_1(x)=2x^2-(3n+4)x+n^2+2n+2$ be a real function in $x$, where $2\leqslant x\leqslant \frac{n}{2}-1$. Consider the derivative of $g_1(x)$, we have
$g'_1(x)=4x-3n-4$. Note that $g'_1(x)\leqslant g'_1(\frac{n}{2}-1)=-n-8<0$, which implies that $g_1(x)$ is a monotonically decreasing function for $2\leqslant x\leqslant \frac{n}{2}-1$. Hence, $g_1(x)\geqslant g_1(\frac{n}{2}-1)=n+8>0$ for $2\leqslant x\leqslant \frac{n}{2}-1$. Furthermore, we have
$\Phi_2(\sqrt{(a-1)(n-a+1)})=-\frac{1}{4}(a-1)(n-2a)g_1(a)<0$, which gives that $\sqrt{(a-1)(n-a+1)}<\rho(G_2)$.

Now we consider that $n$ is odd. If $a=1$, then $G_3=K_{\frac{n-1}{2},\frac{n-1}{2}}\cup K_1$. It is easy to check that
$$
0=\sqrt{(a-1)(n-a+1)}<\frac{n-1}{2}=\rho(G_3),
$$
as desired.

Next we consider $2\leqslant a\leqslant\frac{n-1}{2}$. Recall that $\rho(G_3)$ equals the largest root of $\Phi_3(x)=0$, where $\Phi_3(x)$ is given in \eqref{eq:5.04}.
Substituting $x=\sqrt{(a-1)(n-a+1)}$ into \eqref{eq:5.04} yields
$$
\Phi_3(\sqrt{(a-1)(n-a+1)})=-\frac{1}{4}(a-1)(n-2a+1)(2a^2-(3n+3)a+n^2+n+2).
$$
Let $g_2(x)=2x^2-(3n+3)x+n^2+n+2$ be a real function in $x$, where $2\leqslant x\leqslant \frac{n-1}{2}$. Consider the derivative of $g_2(x)$, we have
$g'_2(x)=4x-3n-3$. Note that $g'_2(x)\leqslant g'_2(\frac{n-1}{2})=-n-5<0$, which implies that $g_2(x)$ is a monotonically decreasing function for $2\leqslant x\leqslant \frac{n-1}{2}$. Therefore, $g_2(x)\geqslant g_2(\frac{n-1}{2})=4>0$ for $2\leqslant x\leqslant \frac{n-1}{2}$. Furthermore, we have
$\Phi_3(\sqrt{(a-1)(n-a+1)})=-\frac{1}{4}(a-1)(n-2a+1)g_2(a)<0$, which deduces that $\sqrt{(a-1)(n-a+1)}<\rho(G_3)$, as required.

By Case 1 and Case 2, we complete the proof.
\end{proof}

\section{\normalsize Some further discussions}\setcounter{equation}{0}
In this paper, we focus on determining the maximum size (resp. the largest spectral radius) of an $n$-vertex graph (resp. an $n$-vertex bipartite graph) without $[a,b]$-factors.
Firstly, we establish a sharp upper bound on the size of a graph with given order forbidding $[a,b]$-factors. Based on this result, we further establish a sharp upper bound on the spectral radius of a graph with given order forbidding $[a,b]$-factors, which proves a stronger version of Cho-Hyun-O-Park's conjecture~\cite[Conjecture 4.4]{CHOP}. In addition, we provide two sharp upper bounds on the size and spectral radius of a bipartite graph with given partite sets forbidding $[a,b]$-factors, respectively. At last, we provide two upper bounds on the size and spectral radius of a bipartite graph with given order forbidding $[a,b]$-factors, respectively. Consequently, we contribute to Problem \ref{pb1} by proving positive results.

In fact, we may view our results by their equivalent forms, in which one may give size condition or spectral condition to guarantee that a graph (or bipartite graph) contains an $[a,b]$-factor. For example, we may reformulate Theorem~\ref{thm1.1} as follows, which presents a size condition to ensure that an $n$-vertex graph contains an $[a,b]$-factor.
\begin{theorem}\label{thm6.1}
Let $a\leqslant b$ be two positive integers, and let $G$ be a graph of order $n$ with $n\geqslant a+1$ and $na\equiv0\pmod{2}$ when $a=b$. If $e(G)\geqslant\binom{n-1}{2}+a-1$, then $G$ contains an $[a,b]$-factor unless one of the following holds:
\begin{wst}
\item[{\rm (i)}]$G\cong K_{a-1}\vee (K_{n-a}\cup K_1)$ or $K_{1,3},$ if $ab=1$ or $ab=2;$
\item[{\rm (ii)}]$G\cong K_{a-1}\vee (K_{n-a}\cup K_1)$ or $K_2\vee \overline{K_3},$ if $a=b=2;$
\item[{\rm (iii)}]$G\cong K_{a-1}\vee (K_{n-a}\cup K_1),$ if $b\geqslant3$.
\end{wst}
\end{theorem}
In what follows, we give an equivalent form of Theorem~\ref{thm1.2}.
\begin{theorem}\label{th6.2}
Let $a\leqslant b$ be two positive integers, and let $G$ be a graph of order $n$ with $n\geqslant a+1$ and $na\equiv0\pmod{2}$ when $a=b$. If $
\rho(G)\geqslant \rho(K_{a-1}\vee (K_{n-a}\cup K_1))
$, then $G$ contains an $[a,b]$-factor unless $G\cong K_{a-1}\vee (K_{n-a}\cup K_1)$.
\end{theorem}

Note that if a graph $G$ contains an $[a,b]$-factor. Then the minimum degree $\delta(G)$ must satisfy $\delta(G)\geqslant a$. However, observing the extremal graphs $G$ in Theorems~\ref{thm6.1}-\ref{th6.2}, we are surprised to find that $\delta(G)=a-1$ for $n>5.$ So it is natural to consider the following interesting and challenging problem.
\begin{problem}\label{pb-2}
  Determine sharp lower bounds on the size or spectral radius of an $n$-vertex graph $G$ with $\delta(G)\geqslant a$ such that $G$ contains an $[a,b]$-factor.
\end{problem}

Recall that let $\mathcal{F}$ be a class of graphs. Denote by $\mathrm{ex}(n,\mathcal{F})$  the maximum number of edges in a connected $\mathcal{F}$-free graph of order $n$. Its corresponding extremal graphs are written as  $\mathrm{Ex}(n,\mathcal{F})$. We also let $\mathrm{Ex}_{sp}(n,\mathcal{F)}$  be the set of connected $\mathcal{F}$-free graphs of order $n$ with maximum spectral radius.

A spanning subgraph of a graph $G$ is its subgraph  whose vertex set is the same as that of $G$. Spanning subgraphs of graphs possessing some given properties, say $H$,  are called $\hat{H}$-\textit{factors}. For example, for a positive integer $k$, a $k$-\textit{factor} of a graph is its spanning subgraph each of whose vertices has a constant degree $k$. So a $1$-factor is a perfect matching. Furthermore a $2$-factor is a set of vertex-disjoint cycles which together cover the vertices of the graph. Here we only consider factors without isolated vertices.


We have the following observation.
\begin{wst}
\item Let $\mathcal{M}$ be the set of all $1$-factors (perfect matchings) on $n$ vertices.\  O \cite{S.O} showed that $\mathrm{ex}(n,\mathcal{M})=\frac{1}{2}n^2-\frac{5}{2}n+5$ for $n\geqslant 10$. Furthermore, according to the results and process of proof in \cite{S.O}, we have  $\mathrm{Ex}_{sp}(n,\mathcal{M})=\mathrm{Ex}(n,\mathcal{M})=\{K_1\vee (K_{n-3}\cup 2K_1)\}$ for large enough $n$.
    
\item Let $\mathcal{H}$ be the set of all connected $2$-factors (Hamilton cycles) on $n$ vertices.\ Ore \cite{OO1961} proved that $\mathrm{ex}(n,\mathcal{H})=\frac{1}{2}n^2-\frac{3}{2}n+2$ and $\mathrm{Ex}(n,\mathcal{H})=\{K_1\vee (K_{n-2}\cup K_1),K_2\vee \overline{K_3}\}$ for $n\geqslant 3$. Fiedler and Nikiforov \cite{FV} showed that $\mathrm{Ex}_{sp}(n,\mathcal{H})=\{K_1\vee (K_{n-2}\cup K_1)\}$. Hence, $\mathrm{Ex}_{sp}(n,\mathcal{H})\subseteq \mathrm{Ex}(n,\mathcal{H})$ for $n\geqslant 3$.
    
\item Let $\mathcal{P}$ be the set of all $\{P_2, P_3, \ldots, \}$-factors on $n$ vertices.\  
    Li and Miao \cite{LM-2021} proved that $\mathrm{ex}(n,\mathcal{P})=\frac{1}{2}n^2-\frac{7}{2}n+9$ for $n\geqslant8$. Moreover, together with the results and process of proof in \cite{LM-2021}, for large enough $n$, we have $\mathrm{Ex}_{sp}(n,\mathcal{P})=\mathrm{Ex}(n,\mathcal{P})=\{K_1\vee (K_{n-4}\cup 3K_1)\}$. 
    
\item Let $\mathcal{K}$ be the set of all $\{K_{1,j}:1\leqslant j\leqslant k\}$-factors on $n$ vertices. \ Miao and Li \cite{ML2023} showed that, for large enough $n$, 
    $\mathrm{ex}(n,\mathcal{K})=\frac{1}{2}n^2-(k+\frac{3}{2})n+\frac{k^2}{2}+\frac{5k}{2}+2$. And according to the results and the process of its proof in \cite{ML2023}, we have $\mathrm{Ex}_{sp}(n,\mathcal{K})= \mathrm{Ex}(n,\mathcal{K})=\{K_1\vee (K_{n-k-2}\cup (k+1)K_1)\}$ for large enough $n$.

\item Let $\mathcal{F}_{a,b}$ be the set of all $[a,b]$-factors on $n$ vertices.\ Let $2\leqslant a\leqslant b$ be two positive integers. By Theorems~\ref{thm1.1} and \ref{thm1.2}, one sees $\mathrm{ex}(n,\mathcal{F}_{a,b})=\frac{1}{2}n^2-\frac{3}{2}n+a$, and  $\mathrm{Ex}_{sp}(n,\mathcal{F}_{a,b})=\mathrm{Ex}(n,\mathcal{F}_{a,b})=\{K_{a-1}\vee (K_{n-a}\cup K_1)\}$ for $n\geqslant \max\{a+1,6\}$ and $na\equiv0\pmod{2}$ when $a=b$.  
    
\item Let $\mathcal{S}=\{F:\text{$F$ is an $n$-vertex graph with $\delta(F)=\delta \geqslant 2$ and $\Delta(F)\leqslant \frac{\sqrt{n}}{40}$}\}$. \ 
    For sufficiently large $n$, Alon and  Yuster \cite{AN} proved that $\mathrm{ex}(n, \mathcal{S})=\frac{1}{2}n^2-\frac{3}{2}n+\delta$, whereas Liu and Ning \cite{LN2023} showed that $\mathrm{Ex}_{sp}(n, \mathcal{S})=\{K_{\delta-1}\vee (K_{n-\delta}\cup K_1)\}$ for sufficiently large $n$. One can verify that $K_{\delta-1}\vee (K_{n-\delta}\cup K_1)\in \mathrm{Ex}(n, \mathcal{S})$. It follows that $\mathrm{Ex}_{sp}(n, \mathcal{S})\subseteq\mathrm{Ex}(n, \mathcal{S})$ for sufficiently large $n$. 
    
\end{wst}
Based on the above observation, we propose the following conjecture.
\begin{conjecture}\label{conj-2}
Let $s$ be a fixed positive integer and $n$ be a sufficiently large integer. Assume that $\mathcal{F}$ is the set of all the $\hat{H}$-factors of $K_n$ for some given property $H$ such that $\mathrm{ex}(n,\mathcal{F})=\frac{1}{2}n^2-(s+\frac{1}{2})n+O(1)$. Then $\mathrm{Ex}_{sp}(n,\mathcal{F)}\subseteq \mathrm{Ex}(n,\mathcal{F})$.
\end{conjecture}


\subsection*{\normalsize Acknowledgements}

We are very grateful to all the referees for their many helpful comments. In particular, their suggestion inspires us to give Problem~\ref{pb-2} in the last section. Shuchao Li received support from the National Natural Science Foundation of China (Grant Nos. 12171190, 11671164) and the Special Fund for Basic Scientific Research of Central Colleges (Grant No. CCNU24JC005).

\end{document}